\documentclass[11pt]{amsart}
\usepackage{amsmath,amsthm,amsfonts,amssymb,latexsym,epsfig}

\newtheorem{theorem}{Theorem}[section]

\newtheorem{lemma}{Lemma}[section]

\newtheorem{cor}{Corollary}[section]

\numberwithin{equation}{section}

\theoremstyle{definition}

\theoremstyle{remark}

\begin{document}
\title{Sums of Powers and Majorization}
\author{Peng Gao}
\address{Department of Computer and Mathematical Sciences,
University of Toronto at Scarborough, 1265 Military Trail, Toronto
Ontario, Canada M1C 1A4} \email{penggao@utsc.utoronto.ca}
\date{January 29, 2007.}


\begin{abstract}
 We study certain sequences involving sums of powers of positive integers and
 in connection with this, we give examples to show that power
 majorization does not imply majorization.
\end{abstract}

\maketitle
\section{Introduction}
\label{sec 1} \setcounter{equation}{0}
  Estimations of sums of powers of positive integers have important applications in
  the study of $l^p$ norms of weighted mean matrices, we leave
  interested readers the recent papers \cite{G} and \cite{Be1}
  for more details in this direction by pointing out here that an
  essential ingredient in \cite{G} is the following lemma of Levin and Ste\v ckin
  \cite[Lemma 1-2, p.18]{L&S}:
\begin{lemma}
\label{lem0}
    For an integer $n \geq 1$,
\begin{eqnarray}
\label{4}
    \sum^n_{i=1}i^r &\geq & \frac {1}{r+1}n(n+1)^r, \hspace{0.1in} 0 \leq r \leq 1, \\
\label{201}
   \sum^n_{i=1}i^r  &\geq & \frac {r}{r+1}\frac
   {n^r(n+1)^r}{(n+1)^r-n^r}, \hspace{0.1in} r \geq 1.
\end{eqnarray}
    Inequality \eqref{201} reverses when $-1 <r \leq 1$.
\end{lemma}
   We note here that in the case $r=0$, the expression on the right-hand side of
   \eqref{201} should be interpreted as the limit of $r \rightarrow
   0$ of the non-zero cases and only the case $r \geq 0$ for \eqref{201} was proved in \cite{L&S}
but one checks easily that the proof extends to the case $r >-1$.

   What we are interested in this paper is to study certain sequences involving sums of powers of
   positive integers. Let
   ${\bf a}= \{a_i\}_{i=1}^\infty$ be an increasing sequence of
positive real numbers.
  We define for any integer $n \geq 1$ and any real number $r$,
\begin{equation*}
  R_{n}(r;{\bf a})= \left(\frac {1}{n} \sum_{i=1}^{n}a_{i}^r\bigg/
\frac {1}{n+1}\sum_{i=1}^{n+1}a_{i}^r\right)^{1/r}, \hspace{0.1in} r
\neq 0; \hspace{0.1in} R_{n}(0;{\bf a})=\frac
{\sqrt[n]{\prod^n_{i=1}a_i}}{\sqrt[n+1]{\prod^{n+1}_{i=1}a_i}}.
\end{equation*}
  For ${\bf a}= \{ i \}_{i=1}^\infty$ being the sequence of positive integers,
  we write $P_n(r)$ for $R_{n}(r; \{ i \}_{i=1}^\infty)$ and we note that for $r > 0$,
  we have the following
\begin{equation}
\label{1}
   \frac {n}{n+1}= \lim_{r \rightarrow +\infty}P_n(r) < P_n(r) < P_n(0)=\frac
   {\sqrt[n]{n!}}{\sqrt[n+1]{(n+1)!}}.
\end{equation}
  The left-hand side inequality above is known as Alzer's inequality
\cite{alz}, and the right-hand side inequality above is known as
Martins' inequality \cite{Mar}.  We refer the readers to
\cite{alz1}, \cite{Xu} and \cite{CGQ} for extensions and refinements
of \eqref{1}. We point out here Alzer considered inequalities
satisfied by $P_n(r)$ for $r<0$ in \cite{alz1} and he showed
\cite[Theorem 2.3]{alz1}:
\begin{equation}
\label{2}
   P_n(0) \leq P_n(r) \leq \lim_{r \rightarrow -\infty}P_n(r)=1.
\end{equation}

  Bennett \cite{Be} proved that for $r \geq 1$,
\begin{equation*}
  P_n(r) \leq P_n(1)=\frac
   {n+1}{n+2}
\end{equation*}
  with the above inequality reversed when $0< r \leq 1$. This inequality
  and inequalities \eqref{1}-\eqref{2} seem to suggest that $P_n(r)$ is a decreasing
  function of $r$. It is the goal of this paper to prove this for $r \leq 1$.
  We will in fact establish this more generally for all
$r$ for $R_{n}(r;{\bf
  a})$  under certain
conditions on the sequence. We will show that the sequence ${\bf a}=
\{ i \}_{i=1}^\infty$ satisfies the condition for $r \leq 1$ and
moreover, $P_n(r) \geq P_n(r')$ for $r'>r, r \leq 1$. The special
case $r=0$ is essentially Martins' inequality.

  Our main tool in this paper is the theory of majorization and
   we recall that for two positive real finite sequences
   ${\bf x}=(x_1, x_2, \ldots, x_n)$ and ${\bf y}=(y_1, y_2,
   \ldots, y_n)$, ${\bf x}$ is said to be
   majorized by ${\bf y}$ if for all convex functions $f$, we have
\begin{equation}
\label{0}
   \sum_{j=1}^{n}f(x_j) \leq \sum_{j=1}^{n}f(y_j).
\end{equation}

    We write ${\bf x} \leq_{maj} {\bf y}$ if this occurs and the
    majorization principle states that if $(x_j)$ and $(y_j)$ are
    decreasing, then ${\bf x} \leq_{maj} {\bf y}$ is equivalent to
\begin{eqnarray*}
\label{10}
   x_1+x_2+\ldots+x_j & \leq & y_1+y_2+\ldots+y_j ~~(1 \leq j \leq
   n-1),
     \\
    x_1+x_2+\ldots+x_n & = & y_1+y_2+\ldots+y_n ~~(n \geq 0).
\end{eqnarray*}
    We refer the reader to \cite[Sect. 1.30]{B&B} for a simple proof of
    this.

    As a weaker notation, we say that ${\bf x}$ is power majorized by ${\bf y}$
  if $\sum_{i=1}^nx_i^p\le\sum_{i=1}^n y_i^p$ for all real $p\not\in[0,1]$ and
  $\sum_{i=1}^nx_i^p\ge\sum y_i^p$ for $p\in[0,1]$.
  We denote power majorization by ${\bf x} \leq_p\,{\bf y}$.
  Clausing \cite{Cl} asked whether ${\bf x} \leq_p\, {\bf y}$
  implies ${\bf x} \leq_{maj} {\bf y}$. Although this is true for $n\le 3$, it is
false in general and counterexamples have been given in \cite{Be0},
\cite{al} and \cite{B&J}. Our study of $P_n(r)$ will also allow us
to give counterexamples to Clausing's question, which will be done
in Section \ref{sec 3}.

\section{The Main Theorem}
\label{sec 2} \setcounter{equation}{0}
\begin{lemma}[{\cite[Theorem 2.4]{M&O&P}}]
\label{lem1}
  If $\alpha_i >0, 1 \leq i \leq n$ and $\beta_1 \geq \beta_2 \geq
  \ldots \geq \beta_n >0$ and $\beta_1 / \alpha_1 \leq \ldots \leq \beta_n /
  \alpha_n$, then $(b_1, \ldots, b_n) \leq_{maj} (a_1, \ldots,
  a_n)$, where $a_i=\alpha_i/\sum^n_{j=1}\alpha_j, b_i=\beta_i/\sum^n_{j=1}\beta_j,
   1 \leq i \leq n$.
\end{lemma}

   We now use this to establish the following
\begin{lemma}
\label{lem2}
  Let $r >0$ and let ${\bf a}= \{a_i\}_{i=1}^\infty$ be an increasing sequence of positive real
  numbers. If for any integer $n \geq 2$, ${\bf a}$ satisfies
\begin{equation}
\label{2.1}
  (n-i)\frac
  {a^r_{n+1-i}}{a^r_{n-i}}+1 \leq (n-i+1)\frac
  {a^r_{n+2-i}}{a^r_{n+1-i}}, \hspace{0.1in} 1 \leq i \leq n-1.
\end{equation}
  Then on writing $\alpha_i=(n+1-i)a^r_{n+2-i}+ia^r_{n+1-i},
  \beta_i=a^r_{n+1-i}, 1 \leq i \leq n$, we have $(c_1, \ldots, c_n) \leq_{maj} (b_1, \ldots,
  b_n)$, where $b_i=\alpha_i/\sum^n_{j=1}\alpha_j, c_i=\beta_i/\sum^n_{j=1}\beta_j,
   1 \leq i \leq n$.

   If ${\bf a}$ satisfies
\begin{equation}
\label{2.1'}
  (i+1)\frac
  {a^r_{i+1}}{a_{i+2}} \leq 1+ i\frac
  {a^r_{i}}{a^r_{i+1}}, \hspace{0.1in} 1 \leq i \leq n-1,
\end{equation}
   Then on writing $\gamma_i=(n+1-i)a^{-r}_{i}+ia^{-r}_{i+1},
  \delta_i=a^{-r}_{i}, 1 \leq i \leq n$, we have $(e_1, \ldots, e_n) \leq_{maj} (d_1, \ldots,
  d_n)$, where $d_i=\gamma_i/\sum^n_{j=1}\gamma_j, e_i=\delta_i/\sum^n_{j=1}\delta_j,
   1 \leq i \leq n$.
\end{lemma}
\begin{proof}
  As the proofs are similar, we will only prove the first assertion of the lemma here.
  It is easy to check that $\alpha_i >0, 1 \leq i \leq n$ and $\beta_1 \geq \beta_2 \geq
  \ldots \geq \beta_n >0$. Hence it suffices to show that $\beta_1 / \alpha_1 \leq \ldots \leq \beta_n /
  \alpha_n$ so that our assertion here follows from Lemma
  \ref{lem1}. Now for $1 \leq i \leq i+1 \leq n$, we have
\begin{equation*}
  \frac {\beta_i}{\alpha_i} = \frac
  {a^r_{n+1-i}}{(n+1-i)a^r_{n+2-i}+ia^r_{n+1-i}}, \hspace{0.1in} \frac {\beta_{i+1}}{\alpha_{i+1}} = \frac
  {a^r_{n-i}}{(n-i)a^r_{n+1-i}+(i+1)a^r_{n-i}}.
\end{equation*}
   From this and our assumption on the sequence, we see that $\beta_i / \alpha_i \leq \beta_{i+1} /
  \alpha_{i+1}$ hold for $1 \leq i \leq i+1 \leq n$ and this
  completes the proof.
\end{proof}
  We note here for any $r>0$, there exists a sequence ${\bf
  a}$ so that either the condition \eqref{2.1} or \eqref{2.1'} is satisfied. For example,
  any positive constant sequence will work. A non-trivial
  example is given in the following
\begin{cor}
\label{cor1}
  Let ${\bf a}= \{ i
  \}_{i=1}^\infty$, then the first assertion of Lemma \ref{lem2} holds for
   $0< r \leq 1$ and the second assertion of Lemma \ref{lem2} holds for any $r>0$.
\end{cor}
\begin{proof}
  As the proofs are similar, we will only prove the first assertion of the corollary here.
  It suffices to check for $0< r \leq 1$,
\begin{equation}
\label{3}
  (n-i)\frac
  {(n+1-i)^r}{(n-i)^r}+1 \leq (n-i+1)\frac
  {(n+2-i)^r}{(n+1-i)^r}, \hspace{0.1in} 1 \leq i \leq n-1.
\end{equation}
   Equivalently, on setting $x=n-i$, it suffices to show $f(x+1)-f(x) \geq
   1$ for $x \geq 1$ with
\begin{equation*}
   f(x)=x\Big ( 1+ \frac 1{x} \Big )^r.
\end{equation*}
   By Cauchy's mean value theorem, we have $f(x+1)-f(x)=f'(\xi)$
   with $x< \xi < x+1$, where
\begin{equation*}
   f'(x)=\Big ( 1+ \frac 1{x} \Big )^{r-1}\Big (1+\frac {1-r}{x} \Big ).
\end{equation*}
   It is easy to see via Taylor expansion that for $0< r \leq 1$,
   $x>0$,
\begin{equation*}
   \Big ( 1+ \frac 1{x} \Big )^{1-r} \leq \Big (1+\frac {1-r}{x} \Big
   ).
\end{equation*}
  We then deduce that $f'(x) \geq 1$ for $x>0$ which completes the proof.
\end{proof}
   Now, we are ready to prove the following
\begin{theorem}
\label{thm1}
  Let $r \neq 0$ be any real number and let ${\bf a}= \{a_i\}_{i=1}^\infty$ be an increasing sequence of positive real
  numbers. For any positive integer $n \geq
  1$,  let ${\bf x}_{n(n+1)}$ to be an $n(n+1)$-tuple,
   formed by repeating $n+1$ times each term of the $n$-tuple:
   $(\frac {a^r_1}{(n+1)\sum^n_{i=1}a^r_i}, \ldots, \frac {a^r_n}{(n+1)\sum^n_{i=1}a^r_i})$ and ${\bf
   y}_{n(n+1)}$ an $n(n+1)$-tuple, formed by repeating $n$ times each term of
   the $(n+1)$-tuple: $(\frac {a^r_1}{n\sum^{n+1}_{i=1}a^r_i}$, $\ldots, \frac
   {a^r_{n+1}}{n\sum^{n+1}_{i=1}a^r_i})$, then if ${\bf a}$ satisfies \eqref{2.1}, ${\bf x}_{n(n+1)} \leq_{maj} {\bf
   y}_{n(n+1)}$ for $r>0$ and if ${\bf a}$ satisfies \eqref{2.1'}, ${\bf x}_{n(n+1)} \leq_{maj} {\bf
   y}_{n(n+1)}$ for $r<0$.
\end{theorem}
\begin{proof}
  As the proofs are similar, we will only prove the case $r>0$ here.
  We note first that here
\begin{eqnarray*}
  && {\bf x}_{n(n+1)}=\Big( x_{i(n+1)+j} \Big )_{ 0 \leq i \leq n-1; 1 \leq j
  \leq n+1}, \hspace{0.1in} x_{i(n+1)+j}=\frac
  {a^r_{n-i}}{(n+1)\sum^n_{i=1}a^r_i}; \\
  && {\bf y}_{n(n+1)}=\Big( y_{in+j} \Big )_{ 0 \leq i \leq n; 1 \leq j
  \leq n}, \hspace{0.1in} y_{in+j}=\frac
  {a^r_{n+1-i}}{n\sum^{n+1}_{i=1}a^r_i}.
\end{eqnarray*}
   It is easy to see that $\sum^{n(n+1)}_{i=1}x_i=\sum^{n(n+1)}_{i=1}y_i$ and
   we need to show that for $1 \leq k \leq n(n+1)-1$,
\begin{equation}
\label{2.2}
   \sum^k_{i=1}x_i \leq \sum^k_{i=1}y_i.
\end{equation}
   It follows from Lemma \ref{lem2} that inequality \eqref{2.2} holds for
   $k=(n+1)i, 1 \leq i \leq n$. Now suppose that there exists a $k_0$ with $(n+1)(j-1) < k_0 < (n+1)j$ for
   some $1 \leq j \leq n$  such that
   inequality \eqref{2.2} holds for all $(n+1)(j-1) < k < k_0$ but fails to hold for $k_0$, then one must have
   $x_{k_0} > y_{k_0}$, but then one checks easily that this implies $x_k > y_k$ for all $k_0 \leq k \leq (n+1)j$,
   which in turn implies that \eqref{2.2} fails to hold for
   $k=(n+1)j$, a contradiction and this means such $k_0$ doesn't
   exist and inequality \eqref{2.2} holds for every $k$ and this
   completes the proof.
\end{proof}

\begin{cor}
\label{cor2}
    Let ${\bf a}= \{a_i\}_{i=1}^\infty$ be an increasing sequence of positive real
  numbers. If it satisfies the relation \eqref{2.1}, then the function $r \mapsto R_{n}(r;{\bf
  a})$ is decreasing for $r \geq 0$. If it satisfies the relation \eqref{2.1'},
  then the function $r \mapsto R_{n}(r;{\bf a})$ is decreasing for $r \leq 0$.
\end{cor}
\begin{proof}
  As the proofs are similar, we will only prove the first assertion
  here. In this case, we may further assume $r>0$ here as the case
  $r=0$ follows from a limiting process.
  Let $r'>r>0$ be fixed and let ${\bf x}_{n(n+1)}$ and ${\bf
  y}_{n(n+1)}$ be two sequences defined as in Theorem \ref{thm1}.
  One then applies \eqref{0} for the convex function $f(u)=u^{r'/r}$
  to conclude that $R_{n}(r;{\bf
  a}) \geq R_{n}(r';{\bf
  a})$. As $r, r'$ are arbitrary, this completes the proof.
\end{proof}
  It now follows from Corollary \ref{cor1} and \ref{cor2} that
\begin{cor}
\label{cor3}
  The function $r \mapsto P_n(r)$ is a decreasing function of $r$ for $r \leq 1$. Moreover,
  $P_n(r) \geq P_n(r')$ for $r'>r, r \leq 1$.
\end{cor}
   We note here the limit case as $r \rightarrow 0$ of Corollary \ref{cor3}
   allows one to obtain $P_n(0) \geq
   P_n(r')$, which is essentially Martins' inequality.

\section{Power Majorization and Majorization}
\label{sec 3} \setcounter{equation}{0}

    Our goal in this section is to give counterexamples to Clausing's
    question mentioned at the end of Section \ref{sec 1}. To achieve
    this, we note here that one can easily deduce from the proof of Corollary
\ref{cor1} that inequality \eqref{3}
   reverses when $r>1$, which means, if we use the notations in
   Lemma \ref{lem2}, that instead of having $(c_1, \ldots, c_n) \leq_{maj} (b_1, \ldots,
  b_n)$, we will have $(b_1, \ldots, b_n) \leq_{maj} (c_1, \ldots,
  c_n)$, where ${\bf b}, {\bf c}$ are constructed as in Lemma
  \ref{lem2} with respect to ${\bf a}=\{i\}_{i=1}^\infty$. This
  further implies that for the sequences ${\bf x}_{n(n+1)}, {\bf
  y}_{n(n+1)}$ constructed as in Theorem \ref{thm1} with respect to ${\bf
  a}=\{i\}_{i=1}^\infty$, we no longer have ${\bf x}_{n(n+1)} \leq_{maj} {\bf
  y}_{n(n+1)}$ for $n \geq 2$. However, if $P_n(r)$ is a decreasing
  function for all $r$, then ${\bf x}_{n(n+1)} \leq_{p} {\bf
  y}_{n(n+1)}$ for all $n \geq 1$, which supplies
  counterexamples to Clausing's question.

  It therefore remains to show that there is at least one $r>1$ such
  that $P_n(r) \geq P_n(r')$ for $r'>r$ and $P_n(r) \leq P_n(r')$ for
  $r' \leq r$. To motivate our approach here, we want to first
  mention a further evidence that supports $P_n(r)$ being a decreasing
  function for all $r$. We note a result of Bennett
  \cite[Theorem 12]{Be1}, which we shall present here in a slightly general form that asserts for
  real numbers $\alpha \geq 1, \beta \geq 1$ and any integer $n \geq
  1$,
\begin{equation*}
   \frac {\Big ( \sum^n_{i=1}i^{\alpha} \Big )\Big ( \sum^n_{i=1}i^{\beta} \Big )}{\sum^n_{i=1}i^{\alpha+
   \beta+1}}
   \geq \frac {\Big ( \sum^{n+1}_{i=1}i^{\alpha} \Big )\Big ( \sum^{n+1}_{i=1}i^{\beta} \Big )}{\sum^{n+1}_{i=1}i^{\alpha+
   \beta+1}}
\end{equation*}
    with the above inequality reversed for $\alpha \leq 1, \beta \leq
    1$. One can easily supply a proof of the above result following
    that of \cite[Theorem 12]{Be1} and we shall leave it to the
    reader. Bennett's result corresponds to the case
    $\alpha=\beta=r$, or explicitly, for $r \geq 1$, $n \geq 1$,
\begin{equation*}
   \frac {\Big ( \sum^n_{i=1}i^r \Big )^2}{\sum^n_{i=1}i^{2r+1}}
   \geq \frac {\Big ( \sum^{n+1}_{i=1}i^r \Big )^2}{\sum^{n+1}_{i=1}i^{2r+1}}
\end{equation*}
   with the above inequality reversed for $r \leq 1$. Now for $r
   \geq 1$, we rewrite the above inequality as
\begin{equation}
\label{3.3}
   \Big ( \frac { \frac 1{n}\sum^n_{i=1}i^r}{ \frac
   1{n+1}\sum^{n+1}_{i=1}i^r} \Big )^2
   \geq \frac {n+1}{n} \frac {\frac {1}{n}\sum^{n}_{i=1}i^{2r+1} }{\frac
   {1}{n+1}\sum^{n+1}_{i=1}i^{2r+1}},
\end{equation}
   and note that it follows from Alzer's inequality (the left-hand side inequality of
    \eqref{1}) that
\begin{equation*}
   \frac {n+1}{n} \frac {\frac {1}{n}\sum^{n}_{i=1}i^{2r+1} }{\frac {1}{n+1}\sum^{n+1}_{i=1}i^{2r+1}}
   \geq \Big ( \frac { \frac 1{n}\sum^n_{i=1}i^{2r+1}}{ \frac
   1{n+1}\sum^{n+1}_{i=1}i^{2r+1}} \Big )^{\frac {2r}{2r+1}}.
\end{equation*}
   We then deduce from this and \eqref{3.3} that $P_n(r) \geq
   P_n(2r+1)$ for $r \geq 1$.

   Bennett's
   result above motivates one to ask in general what can we say
   about the monotonicities of the sequences
\begin{equation*}
  \frac {\Big ( \sum^n_{i=1}i^r \Big
  )^{\alpha}}{\sum^n_{i=1}i^{\alpha(r+1)-1}}, \hspace {0.1in} n=1,
  2, 3, \ldots,
\end{equation*}
   with $\alpha, r$ being any real numbers?

   We now discuss a simple case here which in turn will allow us to achieve our initial goal in this section.
   Before we proceed, we
   note that Bennett used what he called ``the Ratio Principle" to
   obtain his result above. For our purpose in this paper, one can
   regard ``the Ratio Principle" as being equivalent to the
   following lemma in \cite{Xu}:
\begin{lemma}[{\cite[Lemma 2.1]{Xu}}]
\label{lem3}
   Let $\{B_n \}^{\infty}_{n=1}$ and $\{C_n \}^{\infty}_{n=1}$ be strictly increasing positive sequences with
   $B_1/B_2 \leq C_1 / C_2$. If for any integer $n \geq 1$,
\begin{equation*}
  \frac {B_{n+1}-B_n}{B_{n+2}-B_{n+1}} \leq  \frac
  {C_{n+1}-C_n}{C_{n+2}-C_{n+1}}.
\end{equation*}
  Then $B_{n}/B_{n+1} \leq C_{n} / C_{n+1}$ for any integer $n \geq 1$.
\end{lemma}
  We now use this to prove
\begin{theorem}
\label{thm2}
  The sequence
\begin{equation*}
   \frac {\Big ( \sum^n_{i=1}i \Big
  )^{\alpha}}{\sum^n_{i=1}i^{2\alpha-1}}, \hspace {0.1in} n=1,
  2, 3, \ldots,
\end{equation*}
   is increasing for $\alpha \geq 2$ and decreasing for $1 < \alpha < 2$.
\end{theorem}
\begin{proof}
   We need to show now for $n \geq 1$, $\alpha \geq 2$,
\begin{equation*}
  \frac {\Big ( \sum^n_{i=1}i \Big
  )^{\alpha}}{\sum^n_{i=1}i^{2\alpha-1}} \geq \frac {\Big ( \sum^{n+1}_{i=1}i \Big
  )^{\alpha}}{\sum^{n+1}_{i=1}i^{2\alpha-1}},
\end{equation*}
  with the above inequality reversed when $1 < \alpha < 2$.
  We now use Lemma \ref{lem3} to establish this. When $n=1$, this is
  equivalent to show that
\begin{equation*}
  g(\alpha)=1+2^{2\alpha-1}-3^{\alpha}
\end{equation*}
   is greater than or equal to $0$ for $\alpha \geq 2$ and less than
   or equal to $0$ for $1 < \alpha < 2$. It is easy to see that
   $g''(\alpha) \geq 0$ for $\alpha \geq 1$, this combined with the
   observation that $g(1)=g(2)=0$ now establishes our assertion
   above.

   We now prove the theorem for the case $\alpha \geq 2$ and the case $1 < \alpha < 2$
   can be proved similarly. By Lemma \ref{lem3}, it suffices to show for $\alpha \geq 2$,
\begin{equation}
\label{3.4}
   \frac
  {(n+2)^{\alpha}-n^{\alpha}}{(n+1)^{\alpha-1}} \geq  \frac {(n+3)^{\alpha}-(n+1)^{\alpha}}{(n+2)^{\alpha-1}}.
\end{equation}
  We define for $x >0$,
\begin{equation*}
  f(x)=\frac {(x+2)^{\alpha}-x^{\alpha}}{(x+1)^{\alpha-1}},
\end{equation*}
  so that
\begin{equation*}
   f'(x)=\frac {\alpha \Big ((x+2)^{\alpha-1}-x^{\alpha-1} \Big )(x+1)-
   (\alpha-1)\Big ((x+2)^{\alpha}-x^{\alpha} \Big
   )}{(x+1)^{\alpha}}.
\end{equation*}
  By Hadamard's inequality, which asserts that for a continuous convex function $h(x)$ on $[a, b]$,
\begin{equation}
\label{3.6}
   h(\frac {a+b}2) \leq \frac {1}{b-a}\int^b_a h(x)dx \leq \frac
   {h(a)+h(b)}{2},
\end{equation}
   we have for $\alpha \geq 2$,
\begin{equation*}
  \frac {\alpha (x+1)}{\alpha -1} \leq \frac {\alpha }{\alpha
  -1}\frac {1}{(x+2)^{\alpha-1}-x^{\alpha-1}} \int^{(x+2)^{\alpha-1}}_{x^{\alpha-1}}x^{\frac
  {1}{\alpha-1}}dx=\frac
  {(x+2)^{\alpha}-x^{\alpha}}{(x+2)^{\alpha-1}-x^{\alpha-1}}.
\end{equation*}
  This implies that $f(x)$ is a decreasing function for $x>0$, so
  that $f(n) \geq f(n+1)$, which is just \eqref{3.4} and this
  completes the proof.
\end{proof}

   We note here as
\begin{equation*}
  \sum^n_{i=1}i^3= \Big ( \sum^n_{i=1}i \Big
  )^{2},
\end{equation*}
   it follows from this and Theorem \ref{thm2} that
\begin{cor}
\label{cor4}
  The sequence
\begin{equation*}
   \frac {\Big ( \sum^n_{i=1}i^3 \Big
  )^{\alpha}}{\sum^n_{i=1}i^{4\alpha-1}}, \hspace {0.1in} n=1,
  2, 3, \ldots,
\end{equation*}
   is increasing for $\alpha \geq 1$ and decreasing for $1/2 < \alpha < 1$.
\end{cor}

   As one deduces  $P_n(r) \geq P_n(2r+1)$ for $r \geq 1$ from \eqref{3.3},
   it follows from Corollary \ref{cor4} that $P_n(3) \geq P_n(r)$ for $r \geq
   3$ and $P_n(3) \leq P_n(r')$ for $1 < r' < 3$. This combined with
   Corollary \ref{cor3} implies $P_n(3) \leq P_n(r')$ for $r' < 3$.
   Now our discussions above immediately imply that, for example,
\begin{equation*}
   {\bf x}_6= \frac {1}{3(1+2^3)}(1, 1, 1, 2^3, 2^3, 2^3) \leq_{p} \frac {1}{2(1+2^3+3^3)}(1, 1, 2^3, 2^3, 3^3,
   3^3)={\bf y}_6,
\end{equation*}
     and ${\bf x}_6 \leq_{maj} {\bf y}_6$ does not hold, a counterexample to Clausing's question.

\section{Further Discussions}
    We note here that Alzer's inequality (the left-hand side inequality of
    \eqref{1}) can be rewritten as
\begin{equation}
\label{3.1}
    \sum^n_{i=1}i^r \geq  \frac
    {n^{1+r}(n+1)^r}{(n+1)^{1+r}-n^{1+r}}, \hspace{0.1in} r>0.
\end{equation}
    When $0<r \leq 1$, inequality \eqref{3.1} follows from
    \eqref{4}. In fact, one checks easily via the mean value theorem that the
   right-hand side expression in \eqref{4} is greater than or equal to the right-hand
   side expression in \eqref{3.1}. Similarly, when $r \geq 1$, inequality \eqref{3.1}
   follows \eqref{201}.

   Recently, Bennett \cite[Theorem 2]{Be1} has shown that the sequence
\begin{equation*}
  \Big \{ \frac 1{n}\sum^{n}_{i=1}i^r \Big \}^{\infty}_{n=1}
\end{equation*}
   is convex for $r \geq 1$ or $r \leq 0$ and concave for $0 \leq r
   \leq 1$. Equivalently, this is amount to assert that \cite[Theorem 10]{Be1} for $r \geq 1$,
\begin{equation*}
   \sum^{n}_{i=1}i^r \geq \frac {n^r(n+1)^r\Big ((n+2)^r-(n+1)^r \Big
   )}{n^r(n+1)^r-2n^r(n+2)^r+(n+1)^r(n+2)^r},
\end{equation*}
   with the above inequality reversed when $-1 < r \leq 1, r \neq 0$. He then used this to deduce
   that \cite[Corollary
   1]{Be1} for $r \geq 1$,
\begin{equation*}
   \sum^n_{i=1}i^r \geq \frac
    {n^r(n+\frac 1{2})(n+1)^r}{(n+1)^{r+1}-n^{r+1}},
\end{equation*}
   with the above inequality reversed when $-1 < r \leq 1$. We note
   that the above inequality is weaker than inequality \eqref{201} for
   $r > -1$. As an example, we show here for $r \geq 1$,
\begin{equation*}
     \frac {r}{r+1}\frac
    {n^r(n+1)^r}{(n+1)^r-n^r} \geq \frac
    {n^r(n+\frac 1{2})(n+1)^r}{(n+1)^{r+1}-n^{r+1}}.
\end{equation*}
     The above inequality now follows from Hadamard's inequality \eqref{3.6} as
\begin{equation*}
  \frac
    {(n+1)^{r+1}-n^{r+1}}{(n+1)^r-n^r}=
    \frac {r+1}{r}\frac {1}{(n+1)^r-n^r} \int^{(n+1)^r}_{n^r}x^{\frac {1}{r}}dx \geq \frac {r+1}{r}(n+\frac
    1{2}).
\end{equation*}




\begin{thebibliography}{99}
\bibitem{al}
G. D. Allen, Power majorization and majorization of sequences, {\em
Result. Math.}, {\bf 14} (1988), 211-222.
\bibitem{alz}
H. Alzer, On an inequality of H. Minc and L. Sathre, {\em J. Math.
Anal. Appl.}, {\bf 179} (1993), 396-402.
\bibitem{alz1}
H. Alzer, Refinement of an inequality of G. Bennett, {\em  Discrete
Math.}, {\bf 135} (1994), 39-46.
\bibitem{B&B} E. F. Beckenbach and R. Bellman, {\em
Inequalities}, Springer-Verlag, Berlin-G\"ottingen-Heidelberg, 1961.
\bibitem{Be0}
G. Bennett, Majorization versus power majorization, {\em Anal.
Math.}, {\bf 12} (1986), 283-286.
\bibitem{Be}
G. Bennett, Lower bounds for matrices. II., {\em Canad. J. Math.},
{\bf 44} (1992), 54-74.
\bibitem{Be1}
G. Bennett, Sums of powers and the meaning of $l\sp p$, {\em Houston
J. Math.}, {\bf 32} (2006), 801-831.
\bibitem{B&J} G. Bennett and G. Jameson, Monotonic averages of convex
functions, {\em J. Math. Anal. Appl.}, {\bf 252} (2000), 410-430.
\bibitem{CGQ}
T. H. Chan, P. Gao and F. Qi, On a generalization of Martin's
inequality, {\em  Monatsh. Math.}, {\bf 138} (2003), 179-187.
\bibitem{Cl}
A. Clausing, A problem concerning majorization, in {\em General
Inequalities 4} (W. Walter, Ed.), Birkh\"user, Basel, 1984.
\bibitem{G}
P. Gao, {A note on Hardy-type inequalities}, {\it Proc. Amer. Math.
Soc.}, \textbf{133} (2005), 1977-1984.
\bibitem{L&S} V. I. Levin and S. B. Ste\v ckin, Inequalities, {\em
Amer. Math. Soc. Transl. (2)}, {\bf 14} (1960),  1--29.
\bibitem{M&O&P}
A. W. Marshall, I. Olkin and F. Proschan, Monotonicity of ratios of
means and other applications of majorization, {\em Inequalities
(Proc. Sympos. Wright-Patterson Air Force Base, Ohio, 1965)}, pp.
177-190, Academic Press, New York, 1967.
\bibitem{Mar}
J. S. Martins, Arithmetic and geometric means, an application to
Lorentz sequence spaces, {\em Math. Nachr.}, {\bf 139} (1988),
281--288.
\bibitem{Xu}
Z. K. Xu and D. P. Xu, A general form of Alzer's inequality, {\em
Comput. Math. Appl.}, {\bf 44} (2002), 365-373.


\end{thebibliography}
\end{document}